\documentclass[12pt]{amsart}
\usepackage{amsmath,amscd,amssymb,amsfonts}
\newcommand{\h}{\hbox}
\newcommand{\q}{\quad}

\newcommand{\bs}{\par\bigskip}
\newcommand{\ms}{\par\medskip}
\newcommand{\sk}{\par\smallskip}
\newcommand{\bsn}{\par\bigskip\noindent}
\newcommand{\msn}{\par\medskip\noindent}

\newcommand{\mopl}{\hbox{$\bigoplus$}}
\newcommand{\mprod}{\hbox{$\prod$}}
\newcommand{\msum}{\hbox{$\sum$}}
\newcommand{\mcup}{\hbox{$\bigcup$}}
\newcommand{\ot}{\1{\otimes}\1}
\newcommand{\bl}{\bigl}
\newcommand{\br}{\bigl}
\newcommand{\alt}{\widetilde{\alpha}}
\newcommand{\bt}{\widetilde{b}}
\newcommand{\mt}{\widetilde{m}}
\newcommand{\Vt}{\widetilde{V}}
\newcommand{\al}{\alpha}
\newcommand{\la}{\lambda}
\newcommand{\dd}{\partial}
\newcommand{\dt}{\partial_t}
\newcommand{\A}{{\mathcal A}}
\newcommand{\B}{{\mathcal B}}
\newcommand{\BB}{\widetilde{\mathcal B}}
\newcommand{\sD}{\h{$\scriptscriptstyle (*D)$}}
\newcommand{\sDp}{\h{$\scriptscriptstyle (*D')$}}
\newcommand{\BfD}{{\mathcal B}_f^{\sD}}
\newcommand{\OXD}{{\mathcal O}_X^{\sD}}
\newcommand{\OXDp}{{\mathcal O}_{X'}^{\sDp}}
\newcommand{\D}{{\mathcal D}}
\newcommand{\F}{{\mathcal F}}
\newcommand{\G}{{\mathcal G}}
\newcommand{\Hc}{{\mathcal H}}
\newcommand{\I}{{\mathcal I}}
\newcommand{\J}{{\mathcal J}}
\newcommand{\M}{{\mathcal M}}
\newcommand{\OO}{{\mathcal O}}
\newcommand{\Sc}{{\mathcal S}}
\newcommand{\Gt}{\widetilde{\mathcal G}}
\newcommand{\Jt}{\widetilde{\mathcal J}}
\newcommand{\C}{{\mathbb C}}
\newcommand{\N}{{\mathbb N}}
\newcommand{\Q}{{\mathbb Q}}
\newcommand{\Z}{{\mathbb Z}}
\newcommand{\R}{{\mathbb R}}
\newcommand{\can}{{\rm can}}
\newcommand{\mlct}{{\rm mlct}}
\newcommand{\DR}{{\rm DR}}
\newcommand{\Gr}{{\rm Gr}}
\newcommand{\ssb}{\raise.15ex\h{${\scriptscriptstyle\bullet}$}}
\newcommand{\ssc}{\,\raise.2ex\hbox{${\scriptstyle\circ}$}\,}
\newcommand{\into}{\hookrightarrow}
\newcommand{\onto}{\mathop{\rlap{$\to$}\hskip2pt\hbox{$\to$}}}
\newcommand{\simto}{\buildrel{\sim}\over\longrightarrow}
\newcommand{\ges}{\geqslant}
\newcommand{\les}{\leqslant}
\newcommand{\1}{\hskip1pt}
\begin{document}
\title{Hodge ideals and microlocal $V$-filtration}
\author[M. Saito ]{Morihiko Saito}
\begin{abstract} We show that the Hodge ideals in the sense of Mustata and Popa are quite closely related to the induced microlocal $V$-filtration on the structure sheaf, defined by using the microlocalization of the $V$-filtration of Kashiwara and Malgrange. More precisely the former coincide, module the ideal of the divisor, with the part of the latter indexed by positive integers, although they are different without modulo the ideal in general. This coincidence implies that the $j$-log-canonicity in their sense is determined by the microlocal log-canonical threshold of the divisor, which coincides with the maximal root of the reduced (or microlocal) Bernstein-Sato polynomial up to a sign.
\end{abstract}
\maketitle
\centerline{\bf Introduction}
\bsn
Let $D$ be a reduced divisor on a smooth complex algebraic variety $X$ with ${\rm Sing}\,D\ne\emptyset$. For $p\ges 0$, we denote by $\I^{(D,\,p)}\subset\OO_X$ the {\it $p$\1th Hodge ideal} of $D$ in the sense of \cite{MP1}. By definition we have the equalities
$$\I^{(D,\,p)}\bl((p+1)D\br){}=F_p\1\OXD\q\q(\forall\,p\ges 0),$$
where the left-side hand denotes the natural injective image of
$$\I^{(D,\,p)}\otimes_{\OO_X}\OO_X\bl((p+1)D\br)\,\into\OXD:=\OO_X(*D),$$
and $F$ is the {\it Hodge filtration} of mixed Hodge module (see\cite{mhm}) on $\OXD=j_*\OO_U$ with $U:=X\setminus D$ and $j:U\into X$ the natural inclusion. Let $\Vt$ be the {\it microlocal $V$-filtration} on $\OO_X$ defined by using the microlocalization of the $V$-filtration of Kashiwara \cite{Ka2} and Malgrange \cite{Ma2}, see \cite{mic} (and (1.3) below). We have the following.
\msn
{\bf Theorem~1.} {\it Let $\I_D$ be the ideal sheaf of $D$. We have the equalities
$$\I^{(D,p)}=\Vt^{p+1}\OO_X\mod\I_D\q(\forall\,p\ges 0),
\leqno(1)$$
that is, the equalities hold for their images in $\OO_X/\I_D$.}
\ms
This is shown by using the theory of microlocal $V$-filtration, see (2.1) below. The equality holds without mod $\I_D$ in (1) for $p=0,1$, if $D$ is analytic-locally defined by a homogeneous polynomial having an isolated singularity, although this is false for $p=2$, see (2.4) below.
\sk
Let $f$ be a local defining function of $D$ with $b_f(s)$ the Bernstein-Sato polynomial. Set $\bt_f(s):=b_f(s)/(s+1)$. This is called the {\it reduced} (or {\it microlocal}\,) {\it Bernstein-Sato polynomial} of $f$, see (1.3.9) below. Let $-\alt_f$ be the maximal root of $\bt_f(s)$. For instance, if $D$ has a semi-weighted-homogeneous isolated singularity, that is, if $f$ is a $\mu$-constant deformation of a weighted homogeneous polynomial of weights $(w_i)$ having an isolated singularity at the origin, then it is well known that $\alt_f=\msum_i\,w_i$, see (2.5) below.
\sk
One can easily show in general that $\bt_f(s)$ is independent of the choice of $f$ (see for instance \cite[Remark 4.2(i)]{wh}), and the reduced Bernstein-Sato polynomial $\bt_D(s)$ of $D$ can be defined globally. Moreover the maximal root $-\alt_D$ of $\bt_D(s)$ coincides up to a sign with the {\it microlocal log canonical threshold} of $D$, which is by definition the {\it minimal microlocal jumping coefficient} of $D$ (that is defined by using the microlocal $V$-filtration), and is denoted by $\mlct(D)$, see (1.5.5) below. If $\mlct(D)\les 1$, this coincides with the usual {\it log canonical threshold\1} ${\rm lct}(D)$ of $D$, since ${\rm lct}(D)$ coincides up to a sign with the maximal root $-\al_D$ of $b_D(s)$ as is well known (see (1.5.3) below). As a corollary of Theorem~1, we get the following
\msn
{\bf Corollary~1.} {\it We have for any $p\in\N$
$$\I^{(D,\,p)}=\OO_X\iff p\les\mlct(D)-1,
\leqno(3)$$
or equivalently, we have}
$$\min\{p\in\N\mid\I^{(D,\,p)}\ne\OO_X\}=\lfloor\mlct(D)\rfloor.
\leqno(4)$$
\ms
Here $\lfloor\al\rfloor:=\max\{k\in\Z\mid k\les\al\}$.
The implication $\Longleftarrow$ in (3) was shown in \cite[Section 4.5]{bfut}. Note that the equality mod $\I_D$ in Theorem~1 is enough for the proof of (3) in Corollary~1 using (1.3.8) below. In the semi-weighted-homogeneous isolated singularity case explained above, Corollary~1 is closely related to \cite[Theorem 0.9]{mosc}, and gives a generalization of \cite[Theorem D]{MP1} to this case. By Corollary~1 we get the following.
\msn
{\bf Corollary~2.} {\it For $j\in\N$, the pair $(X,D)$ is $j$-{\it log-canonical\1} in the sense of \cite{MP1} $($that is, $\I^{(D,p)}=\OO_X$ for any $p\in[0,j])$ if and only if $j\les\mlct(D)-1$.}
\ms
In the case $j=0$, this is closely related to du Bois singularities, and is already known, see \cite{MP1} and also \cite[Corollary 6.6]{KoSc}, \cite[Remark 4.4]{MSS1}, \cite[Theorem 0.5]{mosc}.
\sk
We can also define the local version of $\mlct(D)$, that is, $\mlct(D,x)$ for $x\in D$. This coincides up to a sign with the maximal root of the reduced local Bernstein-Sato polynomial $\bt_{f,x}(s)=b_{f,x}(s)/(s+1)$. Combining Corollary~1 with \cite{DMST}, we get the following.
\msn
{\bf Corollary~3.} {\it The microlocal log canonical threshold $\mlct(D)$ and the $j$-log-canonicity of $D$ are stable by the restriction to a smooth subvariety of $X$ which is transversal to any strata of a Whitney stratification of $D$. Moreover, $\mlct(D,x)$ is constant on each stratum of a Whitney stratification of $D$.}
\ms
Indeed, let $(\varphi_{f,\la}\OO_X,F)$ be the $\la$-eigensheaf of the underlying filtered $\D_X$-module of the vanishing cycle Hodge module $\varphi_f\Q_{h,X}[d_X-1]$ (up to a shift of filtration) which is defined by using the algebraic partial microlocalization of the direct image of $(\OO_X,F)$ by the graph embedding of $f$, see (1.3.10) below. In particular, we have the isomorphisms
$$\Gr_{\Vt}^{\al}\OO_X=\Gr_F^p(\varphi_{f,\la}\OO_X)\q\q\bl(\la=\exp(-2\pi i\al),\,\,p=-\lfloor 1-\al\rfloor\br).
\leqno(5)$$
Setting $p(f,\la):=\min\bl\{p\in\Z\mid F_p(\varphi_{f,\la}\OO_X)\ne0\br\}$,
we have
$$\aligned\mlct(D)&=\min\bl\{\al\in\Q\mid\Gr_{\Vt}^{\al}\OO_X\ne 0\br\}\\
&=\min\bl\{\al+p(f,\la)\mid\al\in(0,1],\,\la=\exp(-2\pi i\al)\br\}.\endaligned
\leqno(6)$$
The support of $F_p(\varphi_{f,\la}\OO_X)$ is a union of strata of a Whitney stratification of $D$, see (1.4) below. By an argument similar to \cite[lemma 4.2]{DMST}, $F_p(\varphi_{f,\la}\OO_X)$ is compatible with the restriction to a smooth subvariety of $X$ which is transversal to any stratum of the above Whitney stratification, see (1.3.11) below.
\sk
It is also known that the roots of the {\it reduced\,} Bernstein-Sato polynomial are contained in $[\alt_f,d_X-\alt_f]$ with $d_X:=\dim X$ (hence $\mlct(D)=\alt_f\les\lfloor d_X/2\rfloor$), see \cite[Theorem 0.4]{mic} (where $\alt_f$ is denoted by $\al_f$). Combined with Corollary~1, this is quite closely related to \cite[Theorem A (4)]{MP1}. (Indeed, we have $\lfloor n/2\rfloor=\lceil(n-1)/2\rceil$ with $n=d_X$, where $\lceil\al\rceil:=\min\{k\in\Z\mid k\ges\al\}$.)
\sk
This work is partially supported by Kakenhi 15K04816.
\sk
In Section 1 we review some basics of microlocal $V$-filtrations and microlocal multiplier ideals. In Section 2 we prove the main theorem, and describes the induced microlocal $V$-filtration on the structure sheaf in (2.2). In Appendix, we give some remarks related to papers of Mustata and Popa \cite{MP1,MP2}.
\bs\bs
\vbox{\centerline{\bf 1. Microlocal $V$-filtrations and microlocal multiplier ideals}
\bsn
In this section we review some basics of microlocal $V$-filtrations and microlocal multiplier ideals.}
\msn
{\bf 1.1.~$V$-filtration.} Let $f$ be a non-constant function on a smooth complex algebraic variety (or a connected complex manifold) $X$, that is $f\in\Gamma(X,\OO_X)\setminus\C$. Let $i_f:X\into X\times\C$ be the graph embedding by $f$ with $t$ the coordinate of the second factor of $X\times\C$. Set
$$\aligned(\B_f,F):={}&(i_f)_*^{\!\D}(\OO_X,F)\\={}&(\OO_X,F)\ot_{\C}\bl(\C[\dt],F\br),\endaligned$$
where $(i_f)_*^{\!\D}$ denotes the direct image as filtered $\D$-module (up to a shift of filtration), and the set-theoretic direct image is omitted to simplify the notation. The action of $t$ and $\dd_{x_i}$ with $x_i$ local coordinates of $X$ are given by
$$\aligned t\1(g\ot\dt^k)&=fg\ot\dt^k-kg\ot\dt^{k-1},\\ \dd_{x_i}(g\ot\dt^k)&=(\dd_{x_i}g)\ot\dt^k-(\dd_{x_i}f)g\ot\dt^{k+1},\endaligned
\leqno(1.1.1)$$
where $g\in\OO_X$, $k\in\N$ (and the actions of $\dt$ and $\OO_X$ are natural ones).
\sk
More precisely, $g\ot\dt^k$ is an abbreviation of $g\ot\dt^k\delta(t-f)$, and we can identify $\delta(t-f)$ with $f^s$, and $-\dt t$ with $s$ in $\D_X[s]f^s$, see \cite{Ma1}, \cite{Ma2}. We then get the inclusions
$$\OO_X\into\D_X[s]f^s\into\B_f\into\BfD,
\leqno(1.1.2)$$
where $\BfD$ is defined in (1.1.5) below.
\sk
The Hodge filtration $F$ on $\OO_X$ is defined by
$$\Gr_p^F\OO_X=0\q(p\ne 0).$$
To simplify the notation, we do {\it not\1} shift the filtration under the direct image so that
$$F_p\1\B_f=\mopl_{i=0}^p\,\OO_X\ot\dt^i.
\leqno(1.1.3)$$
This is {\it different\1} from the convention in other papers (for instance,\cite{MSS1}, \cite{MSS2}).
\sk
By \cite{mhm}, $\OXD:=\OO_X(*D)$ underlies a mixed Hodge module, denoted by $j_*\Q_{h,X\setminus D}[d_X]$ in this paper, where $j:X\setminus D\into X$ is the natural inclusion and $d_X:=\dim X$. So it has the Hodge filtration $F$. We have the inclusions
$$F_p\1\OXD\subset P_p\1\OXD:=\OO_X\bl((-p-1)D\br)\q\q(\forall\,p\ges 0),
\leqno(1.1.4)$$
by using the Hartogs type property of local sections of $\OO_X$, since the inclusions hold outside the singular locus of $D$ (which has codimension $\ges 2$ in $X$). Note that
$$F_p\1\OXD=P_p\1\OXD=0\q\q(\forall\,p<0).$$
The filtration $P$ is called the {\it pole order filtration}.
\sk
Set
$$\aligned(\BfD,F):={}&(i_f)_*^{\!\D}\bl(\OXD,F\br)\\={}&(\OXD,F)\ot_{\C}\bl(\C[\dt],F\br).\endaligned
\leqno(1.1.5)$$
As in (1.1.3) we have
$$F_p\1\BfD=\mopl_{i=0}^p\,F_{p-i}\OXD\ot\dt^i.
\leqno(1.1.6)$$
\sk
Let $V$ be the filtration of Kashiwara \cite{Ka2} and Malgrange \cite{Ma2} on a regular holonomic $\D_{X\times\C}$-module $\M$ indexed discretely by $\Q$ so that
$$t\1(V^{\al}\M)\subset V^{\al+1}\M,\q\dt(V^{\al}\M)\subset V^{\al-1}\M\q\q(\forall\,\al\in\Q),
\leqno(1.1.7)$$
with $t\1(V^{\al}\M)=V^{\al+1}\M$ ($\forall\,\al>0$), and $\dt t-\al$ is nilpotent on $\Gr_V^{\al}\M$ as in \cite[Definition 3.1.1]{mhp}. Here $V_{\al}=V^{-\al}$, and the $V^{\al}\M$ are locally finitely generated over $\OO_{X\times\C}\langle\dd_{x_i},\dt t\rangle$ with $x_i$ local coordinates of $X$. For instance, $\M$ can be $\B_f$, $\BfD$, $\BfD/\B_f$. Note that the last inclusion of (1.1.2) is strictly compatible with $V$, see \cite[Corollary 3.1.5]{mhp}.
\sk
In the case $\M=\BfD$, we have the bijectivity of
$$t:V^{\al}\BfD\simto V^{\al+1}\BfD\q\q(\forall\,\al\in\Q).
\leqno(1.1.8)$$
This is shown by using the bijectivity of $\Gr_V^{\beta}t:\Gr_V^{\beta}\BfD\simto \Gr_V^{\beta+1}\BfD$ ($\forall\,\beta<0)$ and that of $t$ on $\BfD$ together with an inductive limit argument for $\beta\to-\infty$ in the case $\al\les 0$ (by considering $\BfD/V^{\al}\BfD$). Here the bijectivity of $\Gr_V^{\beta}t$ follows from the nilpotence of the action of $\dt t-\beta$ on $\Gr_V^{\beta}\BfD$. If $\al>0$, it follows from the equality just after (1.1.7).
\sk
Let $j':X\times\C^*\into X\times\C$ be the inclusion. By \cite[Proposition 4.2]{bfut} (using \cite[3.2.2--3]{mhp}), we have the equalities
$$F_p\1\BfD=\msum_{i=0}^p\,\dt^i\bl(V^0\BfD\cap j'_*j'{}^*F_{p-i}\B_f\br)\q\q(\forall\,p\ges 0),
\leqno(1.1.9)$$
using $j'{}^*\BfD=j'{}^*\B_f$.
\msn
{\bf 1.2.~Bernstein-Sato polynomials.} By definition the Bernstein-Sato polynomial $b_f(s)$ coincides with the minimal polynomial of the action of $s=-\dt t$ on
$$\D_X[s]f^s/t(\D_X[s]f^s),$$
where
$$t(\D_X[s]f^s)=\D_X[s]f^{s+1}\subset\BfD,$$
see (1.1.1). Hence $b_f(s)$ can be described by using the filtration $V$ on $\BfD$ together with the filtration $G$ defined by
$$G_k\1\BfD:=t^{-k}(\D_X[s]f^s)\subset\BfD\q\q(k\in\Z).$$
More precisely, if $m_{f,\al}$ denotes the multiplicity of $-\al$ in $b_f(s)$, then, using the bijectivity of (1.1.8), we get
$$m_{f,\al+k}=\min\bl\{i\in\N\mid(s+\al)^i=0\,\,\,\h{on}\,\,\,\Gr^G_k\Gr_V^{\al}\BfD\br\}\q\q(\al>0,\,k\in\N).
\leqno(1.2.1)$$
(Here one may restrict to $\al\in(0,1]$.) Indeed, (1.2.1) is trivial if $k=0$. In general, we have
$$t^k:\Gr^G_k\Gr_V^{\al}\BfD\simto\Gr^G_0\Gr_V^{\al+k}\BfD.
\leqno(1.2.2)$$
\sk
The natural inclusion $\B_f\into\BfD$ induces isomorphisms
$$V^{\al}\B_f\simto V^{\al}\BfD\q\q(\forall\,\al>0).
\leqno(1.2.3)$$
(Indeed, the filtration $V$ is strictly compatible with any morphism of regular holonomic $\D$-modules and $\Gr_V^{\al}(\BfD\!/\B_f)=0$ for $\al\notin\Z_{\les 0}$, see \cite[Corollary 3.1.5 and Lemma 3.1.3]{mhp}.) We then get the filtration $G$ on $\Gr_V^{\al}\B_f$ ($\al\in(0,1]$) such that
$$(\Gr_V^{\al}\B_f,G)\simto(\Gr_V^{\al}\BfD,G)\q\q(\al\in(0,1]).\leqno(1.2.4)$$
\sk
Let $-\al_f$ be the maximal root of $b_f(s)$. By (1.2.1) we get
$$\al_f=\max\bl\{\al\in\Q\mid 1\in V^{\al}\OO_X\br\}{}=\min\bl\{\al\in\Q\mid\Gr_V^{\al}\OO_X\ne 0\br\},
\leqno(1.2.5)$$
with $V$ on $\OO_X$ induced by the inclusions in (1.1.2). Indeed, the first equality follows from
$$1\in V^{\al}\OO_X\iff 1\otimes1\in V^{\al}\B_f\iff \D_X[s]f^s\subset V^{\al}\B_f.$$
\sk
By the negativity of the roots of Bernstein-Sato polynomials \cite{Ka1}, we have $\al_f>0$, or equivalently, $G_0\BfD\subset V^{>0}\BfD$. This means that
$$G_{-1}\Gr_V^{\al}\BfD=0\q(\forall\,\al\in(0,1]).
\leqno(1.2.6)$$
\msn
{\bf 1.3.~Microlocal $V$-filtration and Bernstein-Sato polynomials.} In the above notation, we denote by $\BB_f$ the {\it algebraic partial microlocalization} of $\B_f$ (see \cite{mic}), that is,
$$(\BB_f,F)=\OO_X\ot_{\C}(\C[\dd_{t},\dd_{t}^{-1}],F),\q\h{so that}$$
$$\Gr_p^F\BB_f=\OO_X\ot\dd_{t}^{\,p}\q(p\in\Z).
\leqno(1.3.1)$$
We have the {\it microlocal $V$-filtration} $V$ on $\BB_f$ along $t=0$, which is defined by modifying the $V$-filtration of Kashiwara \cite{Ka2} and Malgrange \cite{Ma2} on $\B_f$ as follows (see \cite[2.1.3]{mic}): 
$$V^{\al}\BB_f:=\begin{cases}V^{\al}\B_f\oplus(\OO_X\ot_{\C}\1\C[\dt^{-1}]\dt^{-1})&\h{if}\,\,\,\,\al\les 1,\\ \dt^{-j}V^{\al-j}\BB_f&\h{if}\,\,\,\,\al>1,\,\,\al-j\in(0,1].\end{cases}$$
This is an exhaustive decreasing filtration indexed discretely by $\Q$. We have the natural inclusion
$$\can:\B_f\into\BB_f,$$
preserving the filtrations $F,V$, and inducing the filtered isomorphisms (see \cite[2.1.4]{mic}):
$$\can:\Gr_V^{\al}(\B_f,F)\simto\Gr_V^{\al}(\BB_f,F)\q\q(\forall\,\al<1),
\leqno(1.3.2)$$
together with the bifiltered isomorphisms (see \cite[Lemma 2.2]{mic}):
$$\dt^{\,k}:(\BB_f;F,V)\simto(\BB_f;F[-k],V[-k])\q\q(\forall\,k\in\Z).
\leqno(1.3.3)$$
\sk
Define the filtration $G$ on $\BB_f$ by
$$G_k\BB_f:=\D_X[s,\dt^{-1}](1\ot\dt^k)\subset\BB_f\q\q(k\in\Z),$$
where $\D_X[s,\dt^{-1}]$ denotes the subalgebra of $\D_X[t]\langle\dt,\dt^{-1}\rangle$ generated by $\D_X[s]$ and $\dt^{-1}$ (with $s:=-\dt t$). The {\it microlocal Bernstein-Sato polynomial} $\bt_f(s)$ is defined to be the minimal polynomial of the action of $s$ on $\Gr^G_0\BB_f$. Its maximal root $-\alt_f$ is given up to a sign by
$$\aligned\alt_f={}&\min\bl\{\al\in\Q\mid\Gr_V^{\al}\Gr^G_0\BB_f\ne 0\br\}\\={}&\max\bl\{\al\in\Q\mid\Gr_V^{\beta}\Gr^G_0\BB_f=0\,\,\,(\forall\,\beta<\al)\br\}\\ \buildrel{\scriptscriptstyle(*)}\over{=}{}&\max\bl\{\al\in\Q\mid G_0\BB_f\subset V^{\al}\BB_f\br\}\\={}&\max\bl\{\al\in\Q\mid 1\ot 1\in V^{\al}\BB_f\br\}.\endaligned
\leqno(1.3.4)$$
For the proof of the third equality $(*)$, we use Zassenhaus lemma (see \cite[1.2.1]{De2}) asserting
$$\Gr_V^{\beta}\Gr^G_0\BB_f=V^{\beta}G_0\BB_f/(V^{>\beta}G_0\BB_f+V^{\beta}G_{-1}\BB_f).
\leqno(1.3.5)$$
We have the inclusion $G_0\BB_f\subset V^{\beta}\BB_f$ for $\beta\ll 0$. In the case $G_0\BB_f\subset V^{\beta}\BB_f$, we have
$$\Gr_V^{\beta}\Gr^G_0\BB_f=G_0\BB_f/V^{>\beta}G_0\BB_f,$$
(since $G_{-1}\BB_f=\dt^{-1}G_0\BB_f\subset\dt^{-1}V^{\beta}\BB_f\subset V^{>\beta}\BB_f$ in this case), and hence
$$\Gr_V^{\beta}\Gr^G_0\BB_f=0\iff G_0\BB_f\subset V^{>\beta}\BB_f.
\leqno(1.3.6)$$
So the third equality $(*)$ follows. (The proofs of the other equalities are easy.)
\sk
Let $\mt_{f,\al}$ be the multiplicity of $-\al$ in $\bt_f(s)$. By arguments similar to the proofs of (1.2.1) and (1.2.5), we have
$$\mt_{f,\al+k}=\min\bl\{i\in\N\mid s^i=0\,\,\,\h{on}\,\,\,\Gr_V^{\al}\Gr^G_k\BB_f\br\}\q\q(\al\in\Q,\,k\in\Z),
\leqno(1.3.7)$$
together with the equality
$$\aligned\alt_f=\al':={}&\max\bl\{\al\in\Q\mid 1\in\Vt^{\al}\OO_X\br\}\\ ={}&\min\bl\{\al\in\Q\mid\Gr_{\Vt}^{\al}\OO_X\ne 0\br\},\endaligned
\leqno(1.3.8)$$
where the filtration $\Vt$ on $\OO_X$ is induced from the microlocal $V$-filtration on $\BB_f$ by using the isomorphism (1.3.1) for $p=0$. Indeed, the inequality $\alt_f\les\al'$ holds by definition, and we have by an argument similar to (1.3.5--6) (with $G$ replaced by $F$)
$$\Gr_V^{\alt_f}\Gr^F_0\BB_f=F_0\BB_f/V^{>\alt_f}F_0\BB_f\ne0,$$
since it contains $1\ot 1\notin V^{>\alt_f}\BB_f$.
\sk
Note also that we have by \cite[Theorem 0.3]{mic} the equality
$$\bt_f(s)=b_f(s)/(s+1),
\leqno(1.3.9)$$
that is, the {\it microlocal} Bernstein-Sato polynomial coincides with the {\it reduced} one.
\sk
For $\la\in\C^*$, set
$$\varphi_{f,\la}(\OO_X,F):=\Gr_V^{\al}(\BB_f,F)\q\bl(\al\in(0,1],\,\la=\exp(-2\pi i\al)\br).
\leqno(1.3.10)$$
Then (5) in the introduction follows from (1.3.3). The direct sum of the above filtered $\D_X$-modules underlies the vanishing cycle mixed Hodge module (up to a shift of filtration) which is denoted by $\varphi_f\1\Q_{h,X}[d_X-1]$ in this note, where $d_X:=\dim X$. Its underlying $\Q$-complex is the vanishing cycle complex (see \cite{De3}) which is denoted by $\varphi_f\1\Q_X[d_X-1]$ in this note.
\sk
By an argument similar to \cite[lemma 4.2]{DMST}, we have the canonical isomorphisms
$$F_p(\varphi_{f_Z,\la}\OO_Z)=i_Z^*F_p(\varphi_{f,\la}\OO_X)\q\q(p\in\Z),
\leqno(1.3.11)$$
if $Z\subset X$ is a smooth subvariety transversal to any stratum of a Whitney stratification of $D$, where $i_Z:Z\into X$ is the natural inclusion, and $f_Z:=f|_Z$.
\msn
{\bf 1.4.~Support of Hodge filtrations.} Let $\M$ be a regular holonomic $\D_X$-module, and $\F\subset\M$ be a coherent $\OO_X$-submodule. It is well known that there is a Whitney stratification $\Sc$ of ${\rm Supp}\,\M$ such that
$${\rm Ch}\,\M\subset\mcup_{S\in\Sc}\,T_S^*X,
\leqno(1.4.1)$$
where the left-hand side is the characteristic variety of $\M$, and $T_S^*X$ is the conormal bundle of a stratum $S\in\Sc$ in $X$.
Then $Z:={\rm Supp}\,\F$ must be a union of strata of $\Sc$.
\sk
Indeed, let $\M'\subset\M$ be the $\D_X$-submodule generated by $\F$. Then $Z={\rm Supp}\,\M'$. By the property of characteristic varieties, (1.4.1) holds with $\M$ replaced by $\M'$. Let $S\in\Sc$ be a locally maximal dimensional stratum with $Z\cap S\ne\emptyset$ (assumed connected). Then we must have $Z\cap S=S$, since (1.4.1) cannot hold for $\M'$ unless $Z\cap S=S$ (by considering at smooth points of $Z\cap S$). So the assertion follows, since $Z$ is a closed subset.
\msn
{\bf 1.5.~Multiplier ideals.} In the notation of (1.1--2), let $\J(\al D)\subset\OO_X$ be the {\it multiplier ideals} of $D:=f^{-1}(0)\subset X$ with coefficients $\al$ in $\R$. These can be defined by local integrability of $|g|/|f|^{2\al}$ for $g\in\OO_X$, see \cite{La}, \cite{Na}. They form a weakly decreasing family of coherent ideal sheaves of $\OO_X$, and their graded quotients $\G(\al D)$ {\it can be defined} by
$$\G(\al D):=\J((\al-\varepsilon)D)/\J(\al D)\q(0<\forall\,\varepsilon\ll 1),$$
with $\J((\al+\varepsilon)D)=\J(\al D)$ for $0<\forall\,\varepsilon\ll 1$ (where the range of $\varepsilon$ depends on $\al$). Using a resolution of $f$, one can show that
$${\rm JC}(D):=\bl\{\al\in\R\mid\G(\al D)\ne 0\br\}{}\subset\Q,
\leqno(1.5.1)$$
and the members of ${\rm JC}(D)$ are called the {\it jumping coefficients} of $D$, see \cite{La}. We then assume that the $\al$ are rational numbers when we consider $\J(\al D)$, $\G(\al D)$.
\sk
By \cite[Theorem 0.1]{BS} we have
$$\aligned\J(\al D)&=V^{\al}\OO_X&\h{if}\,\,\,\al\notin{\rm JC}(D),\\ \G(\al D)&=\Gr_V^{\al}\OO_X=V^{\al}\OO_X/\J(\al D)&\h{if}\,\,\,\al\in{\rm JC}(D),\endaligned
\leqno(1.5.2)$$
where the filtration $V$ on $\OO_X$ is induced by the filtration $V$ of Kashiwara and Malgrange on $\B_f$ via the inclusion (1.1.2). By (1.2.5) and (1.5.2), we get in the notation of the introduction
$$\al_D={\rm lct}(D)\,\bl(:=\min{\rm JC}(D)\br).
\leqno(1.5.3)$$
\sk
We define the {\it microlocal multiplier ideals} $\Jt(\al D)$ and their graded quotients $\Gt(\al D)$ so that
$$\aligned\Jt(\al D)&=\Vt^{\al}\OO_X&\h{if}\,\,\,\al\notin\widetilde{\rm JC}(D),\\ \Gt(\al D)&=\Gr_{\Vt}^{\al}\OO_X=\Vt^{\al}\OO_X/\Jt(\al D)&\h{if}\,\,\,\al\in\widetilde{\rm JC}(D).\endaligned
\leqno(1.5.4)$$
where $\widetilde{\rm JC}(D):=\{\al\in\Q\mid\Gt(\al D)\ne0\}$ (called the set of {\it microlocal jumping coefficients} of $D$), and the filtration $\Vt$ on $\OO_X$ is as in (1.3.8), see also \cite{MSS2}. Note that we have for $0<\varepsilon\ll 1$ (depending on $\al$)
$$\Gt(\al D)=\Jt((\al-\varepsilon)D)/\Jt(\al D),\q\Jt((\al+\varepsilon)D)=\Jt(\al D).$$
By (1.3.8) and (1.5.4) we get in the notation of the introduction
$$\alt_D=\mlct(D)\,\bl(:=\min\widetilde{\rm JC}(D)\br).
\leqno(1.5.5)$$
\bs\bs
\vbox{\centerline{\bf 2. Proof of the main theorem}
\bsn
In this section we prove the main theorem, and describes the induced microlocal $V$-filtration on the structure sheaf in (2.2).}
\msn
{\bf 2.1.~Proof of Theorem~1.} Since the assertion is local, we may assume that $D$ is defined by a function $f$. In the notation of (1.1), set
$$v^{(p)}:=\sum_{k=0}^p\frac{1}{\h{$f^{p+1-k}$}}\ot\frac{\dt^{\1 k}}{k!}\in\BfD\q\q(p\in\N).$$
By (1.1.1) we have
$$t\1 v^{(p)}=\sum_{k=0}^p\frac{1}{\h{$f^{p-k}$}}\ot\frac{\dt^{\1 k}}{k!}\,-\,\sum_{k=1}^p\frac{1}{\h{$f^{p+1-k}$}}\ot\frac{\dt^{\1 k-1}}{(k-1)!}=1\ot\frac{\dt^{\1 p}}{p!}\in\B_f\q(p\in\N).
\leqno(2.1.1)$$
\sk
Let $g_0\in\I^{(D,p)}$, that is, $g_0/f^{p+1}\in F_p\OXD$. By the direct sum decomposition (1.1.6) together with the formula (1.1.9), there is
$$w=\sum_{k=0}^p\dfrac{g_k}{\h{$f^{p+1-k}$}}\ot\dfrac{\dt^k}{k!}\in V^0P_p\1\BfD\q\h{with}\q g_k\in\OO_X\,\,(k\in[1,p]),$$
where $P_p\1\BfD$ is defined in the same was as in (1.1.6) with $F$ replaced by $P$. By the first inclusion of (1.1.7) for $\al=0$ together with the isomorphism (1.2.3) for $\al=1$, we have
$$tw=\msum_{k=0}^p\,h_k\ot\dt^k/k!\in V^1\1\BfD=V^1\1\B_f\q\h{with}\q h_k\in\OO_X\,\,\bl(k\in[0,p]\br).
\leqno(2.1.2)$$
By (2.1.1) together with the injectivity of the action of $t$ on $\BfD$, we then get
$$w=\msum_{k=0}^p\,h_k\1 v^{(k)}\q\h{in}\,\,\,\BfD.$$
In particular
$$g_0=\msum_{k=0}^p\,f^{p-k}h_k.
\leqno(2.1.3)$$
By (2.1.2) together with the filtered isomorphism (1.3.3) for $k=d$, we have
$$h_p\ot 1/p!=\Gr_F^0(\dt^{-p}tw)\in V^{p+1}\1\Gr^F_0\BB_f=\Vt^{p+1}\OO_X.
\leqno(2.1.4)$$
By (2.1.3--4) we thus get
$$g_0\in\Vt^{p+1}\OO_X+\I_D.$$
\sk
Conversely, let $h_p\in\Vt^{p+1}\OO_X$. There are $h_k\in\OO_X$ ($k\in[0,p-1]$) with
$$w':=\msum_{k=0}^p\,h_k\ot\dt^k/k!\in V^1\1\B_f,
\leqno(2.1.5)$$
by the definition of $V^1\BB_f$ written just after (1.3.1). (Note that $F_{-1}\BB_f\subset V^1\BB_f$.) By (2.1.1) we have
$$w'=tw\q\h{with}\q w:=\msum_{k=0}^p\,h_k\1 v^{(k)}\in V^0F_p\1\BfD,
\leqno(2.1.6)$$
where the formula (1.1.9) is used to show the last inclusion. Write
$$w=\sum_{k=0}^p\dfrac{g_k}{\h{$f^{p+1-k}$}}\ot\dfrac{\dt^k}{k!}\q\h{with}\q g_k\in\OO_X\,\,(k\in[0,p]).$$
By the definition of $v^{(k)}$, we have
$$g_0=\msum_{k=0}^p\,f^{p-k}h_k.$$
By (2.1.6) together with the direct sum decomposition (1.1.6), we get
$$g_0/f^{p+1}\in F_p\OXD,\q\h{that is,}\q g_0\in\I^{(D,p)}.$$
These imply that
$$h_p\in\I^{(D,p)}+\I_D.$$
This finishes the proof of Theorem~1.
\msn
{\bf 2.2.~Calculation of $\Vt$ in the weighted homogeneous case.} Let $f$ be a weighted homogeneous polynomial of weights $(w_i)$ having an isolated singularity at the origin; that is, $f$ is a linear combination of monomials $\mprod_i\,x_i^{m_i}$ with $\sum_im_iw_i=1$, where $w_i>0$ ($\forall\,i$). Let $(v_j)_{j\in[1,\mu]}$ be a monomial basis of the Milnor ring $\OO_{X^{\rm an},0}/(\dd f)$, where $(\dd f)$ denotes the Jacobian ideal generated by the partial derivatives $f_i:=\dd f/\dd x_i$, and $\mu$ is the Milnor number of $f$. Set
$$A:=\C\{x_1,\dots,x_n\}\,\bl(=\OO_{X^{\rm an},0}\br),\q B:=\C\{y_1,\dots,y_n\}.$$
Then we have the direct sum decomposition as $B$-modules
$$A=\mopl_{j=1}^{\mu}\,B\1v_j,
\leqno(2.2.1)$$
that is, $A$ is a $B$-module freely generated by the $v_j$ under the $\C$-algebra morphism $B\into A$ defined by $y_i=f_i$, which corresponds to the finite flat surjective morphism
$$(f_1,\dots,f_n):(\C^n,0)\to(\C^n,0).
\leqno(2.2.2)$$
Set
$$\aligned&y^{\nu}:=\mprod_i\,f_i^{\nu_i}\in A,\q|\nu|:=\msum_i\,\nu_i\q\bl(\nu=(\nu_1,\dots,\nu_n)\in\N^n\br),\\&\Lambda(\al):=\{(j,\nu)\in[1,\mu]\times\N^n\mid\al(v_j)+|\nu|\ges\al\}\q(\al\in\Q),\\&\h{with}\q\al(v):=\msum_{i=1}^n\,(m_i+1)w_i\,\,\,\h{for a monomial}\,\,\,v=\mprod_i\,x_i^{m_i}.\endaligned
\leqno(2.2.3)$$
We {\it assume} that $f$ contains monomials of type $x_i^{a_i}$ (with $a_i=1/w_i$) for any $i\in[1,n]$. In the case of homogeneous polynomials (that is, $w_i=1/d$ for any $i$), this assumption is satisfied by replacing coordinates if necessary.)
\msn
{\bf Proposition.} {\it With the above notation and assumptions, we have}
$$\Vt^{\al}\!A=\msum_{(j,\nu)\in\Lambda(\al)}\,A\1y^{\nu}v_j\q(\forall\,\al\in\Q).
\leqno(2.2.4)$$
\msn
{\bf Remark.} This is a generalization of \cite[Example 2.6]{MSS2} where the assertion is proved in the case $f=\sum_ix^{a_i}$ using the Thom-Sebastiani type theorem for microlocal $V$-filtrations.
\msn
{\bf 2.3.~Proof of Proposition in (2.2).} By \cite[Proposition 3.2]{mic}, each $V^{\al}\BB_f$ is generated over $\D_{X,0}\langle\dt^{-1},s\rangle$ by $v\ot\dt^{-k}$ with $v$ a monomial and $k\in\Z$ satisfying
$$\al(v)+k\ges\al.$$
We then get the inclusion $\supset$ in (2.2.4) by considering
$$\mopl_i\,\dd_{x_i}^{\nu_i}\bl(v\ot\dt^{-|\nu|}\br),$$
for $v$ satisfying the above inequality. So the assertion is reduced to the equality
$$\dim\Gr_{\Vt}^{\al}A=\dim\Gr_{V'}^{\al}A\q(\forall\,\al\in\Q),
\leqno(2.3.1)$$
where $V^{\prime\,\al}A$ is defined by the right-hand side of (2.2.4).
\sk
Let $E$ be the filtration on $A$ defined by
$$E^k:=\msum_{|\nu|=k}\,A\1y^{\nu}\q(k\in\Z).$$
Let $V''$ be the filtration on $A$ such that $V''{}^{\1\al}\!A$ is generated over $A$ by monomials $v$ with $\al(v)\ges\al$.
By the definition of $V'$, we see that the filtration $V'$ on $\Gr_E^k$ is given by the strict surjection
$$\msum_{|\nu|=k}\,y^{\nu}:\mopl_{|\nu|=k}(A,V''[-k])\onto\Gr_E^k(A,V'),
\leqno(2.3.2)$$
since the multiplication by $f_i:=\dd_{x_i}f$ induces the injective strict morphism
$$f_i:(A,V'')\into(A,V''[1-w_i])\q(\forall\,i),$$
where $V''[\beta]^{\al}=V''{}^{\al+\beta}$. (Indeed, $f_i$ preserves, up to the shift by $1-w_i$, the grading of $A$ defined by $\deg x=w_i$.)
\sk
By (2.2.1) the kernel of the surjection in (2.3.2) is given by the direct sum of the Jacobian ideals $(\dd f)\subset A$ indexed by $\nu\in\N^n$ with $|\nu|=k$. We then get
$$\dim\Gr_{V'}^{\al}\Gr_E^kA=\msum_{|\nu|=k}\,n_{f,\al-k}=\tbinom{n+k-1}{n-1}\,n_{f,\al-k},$$
with
$$\aligned n_{f,\al}:=\#\{j\mid\al(v_j)=\al\}=\dim\Gr_F^pH^{n-1}(F_{\!f},\C)_{\la}&\\ \bl(\la=\exp(-2\pi i\al),\,p=\lfloor n-\al\rfloor\br),\,\,\,\q\q\q&\endaligned
\leqno(2.3.3)$$
where $F_{\!f}$ denotes the Milnor fiber of $f$, $H^{n-1}(F_{\!f},\C)_{\la}$ is the $\lambda$-eigenspace of the Milnor cohomology, and $F$ is the Hodge filtration, see \cite{St}. (Note that the second equality of (2.3.3) follows from \cite{ScSt,Va1}.) So $n_{f,\al}$ is the multiplicity of the {\it Steenbrink spectrum}
$${\rm Sp}(f):=\msum_{j=1}^{\mu}\,t^{\1\al(v_j)}=\msum_{\al\in(0,n)}\,n_{f,\al}\1 t^{\1\al}.$$
Note that the mixed Hodge structure on the Milnor cohomology $H^{n-1}(F_{\!f},\Q)$ is identified with the vanishing cycle Hodge module $\varphi_f\Q_{h,X}[n-1]$ supported at 0 (see \cite{mhm}) by using a projective compactification of the morphism $f$ as in \cite{Br}. We then get
$$\dim\Gr_{V'}^{\al}A=\msum_{\nu\in\N^n}\,n_{f,\al-|\nu|}=\msum_{k\ges 0}\,\tbinom{n+k-1}{n-1}\,n_{f,\al-k},
\leqno(2.3.4)$$
\sk
On the other hand we have by the definition $\Vt$ and using the bifiltered isomorphism (1.3.3) together with the filtered isomorphism (1.3.10)
$$\aligned&\Gr_{\Vt}^{\al+p}A=\Gr_V^{\al+p}\Gr^F_0\BB_{f,0}=\Gr_V^{\al}\Gr^F_{k+p}\BB_{f,0}=\Gr^F_p\varphi_{f,\la}\OO_{X,0}\\ &\q\q\h{for}\,\,\,\,\al\in(0,1],\,\,p\in\N,\,\,\la=\exp(-2\pi i\al).\endaligned
\leqno(2.3.5)$$
We have moreover the isomorphisms as filtered $\D_{X,0}$-modules for $\al\in(0,1]$, $\la=\exp(-2\pi i\al)$
$$(\varphi_{f,\la}\OO_{X,0},F)=\mopl_{q\in[0,n-1]}\,\mopl_{\al(v_j)=\al+q}\,\bl(\C[\dd_1,\dots,\dd_n](v_j\ot 1),F[q]\br),
\leqno(2.3.6)$$
where $\dd_i:=\dd_{x_i}$, the filtration $F$ on $\C[\dd_1,\dots,\dd_n]\subset\D_{X,0}$ is induced by $F$ on $\D_{X,0}$ (that is, by the order of $\dd_1,\dots,\dd_n$).
This follows from the compatibility of the Hodge structure on the Milnor cohomology $H^{n-1}(F_{\!f},\Q)$ with the vanishing cycle Hodge module $\varphi_f\Q_{h,X}[n-1]$ as is explained before (2.3.4). Note that there is a {\it shift} of the filtration $F$ by $n-1$ between $(H^{n-1}(F_{\!f},\C),F)$ and $(\varphi_{f,\la}\OO_{X,0},F)$, where $F_p=F^{-p}$. (Indeed, compare (2.3.3) with (5) in the introduction). The polynomial ring $\C[\dd_1,\dots,\dd_n]$ appears in (2.3.6) by the definition of the direct image of filtered $\D$-modules under the inclusion $i_0:\{0\}\into X$. (Here the filtration $F$ is not shifted under the direct image by $i_0$.)
\sk
Comparing (2.3.3--4) with (2.3.5--6), we then get the desired equality (2.3.1). This finishes the proof of Proposition in (2.2).
\msn
{\bf 2.4.~Calculation of Hodge ideals in the homogeneous case.} In the notation and assumption of (2.2), we have by \cite[Theorem 0.7]{mosc})
$$F_p\1\OO_{X^{\rm an},0}^{\sD}=\msum_{|\nu|\les p}\,\dd^{\nu}_x\bl(A^{\ges p+1-|\nu|}/f^{p+1-|\nu|}\br)\q(\forall\,p\in\N),
\leqno(2.4.1)$$
where $\dd_x^{\nu}:=\prod_i\dd_i^{\nu_i}$, and $A^{\ges k}\subset A$ is the ideal generated by monomials $v$ with $\al(v)\ges k$. However, it is not easy to calculate further this in general.
\sk
Assume $w_i=1/d\,\,\,(\forall\,i)$ for some positive integer $d$ so that $\al(v)=(\deg v+n)/d$. Then we can show the following by using an argument is similar to the proof of \cite[Lemma 1.5]{DSW}:
$$(\I^{(D,p)})^{\rm an}_0=\Vt^{p+1}\!A\,(=\msum_{(j,\nu)\in\Lambda(p+1)}\,A\1y^{\nu}v_j\br)\q\h{for}\,\,\,\,p=0,1,
\leqno(2.4.2)$$
Indeed, the assertion for $p=0$ follows from (2.4.1) immediately. For $p=1$, the assertion is easily verified if $d\les n$. In the case $d>n$ (that is, $\alt_f=n/d<1$), we have the following by (2.4.1) for {\it any monomial} $v$ with $\al(v)=1-\tfrac{1}{d}$ (that is, $\al(x_iv)=1$):
$$\dd_j(x_iv/f)=(\dd_jx_iv)/f-f_jx_iv/f^2\in F_1\1\OO_{X^{\rm an},0}^{\sD}\q(\forall\,i,j\in[1,n]).
\leqno(2.4.3)$$
Adding this for $i=j\in[1,n]$, and using $\tfrac{1}{d}\sum_if_ix_iv=fv$, we get
$$\tfrac{1}{d}\,\msum_i\,\dd_i(x_iv/f)=\bl(\al(v)-1\br)v/f\in F_1\1\OO_{X^{\rm an},0}^{\sD}\,\,\,\bl(\h{that is,}\,\,\,fv\in(\I^{(D,1)})^{\rm an}_0\br).
\leqno(2.4.4)$$
The condition $\alt_f=n/d<1$ implies that any monomial $\widetilde{v}$ with $\al(\widetilde{v})=1$ is written as $x_{i_0}v$ for some monomial $v$ and $i_0\in[1,n]$. By (2.4.3) with $i=i_0$, we then get the inclusion
$$(\I^{(D,1)})^{\rm an}_0\supset\Vt^{2}\!A,
\leqno(2.4.5)$$
since we have $f(\dd_jx_{i_0}v)\in(\I^{(D,1)})^{\rm an}_0$ by (2.4.4) applied to $v':=\dd_jx_{i_0}v$ instead of $v$. (Note that $v'$ is also a monomial with $\al(v')=1-\tfrac{1}{d}$.)
\sk
Similarly the opposite inclusion $\subset$ holds in (2.4.5), since $fv\in\Vt^2A$ for {\it any} monomial $v$ with $\al(v)=1-\tfrac{1}{d}$ (by using $\tfrac{1}{d}\sum_if_ix_iv=fv$). So (2.4.2) for $p=1$ is proved.
\msn
{\bf Remarks.} (i) The equality (2.4.2) holds for any $p\in\N$ in the case $d=2$ (that is, the singularity is an ordinary double point) by comparing (2.2.4) with \cite[Lemma 1.5]{DSW}.
\ms
(ii) The equality (2.4.2) is false for $p=2$ in the case $f=x^3+y^3+z^3$ with $n=d=3$. Indeed, we have $\al(1)=1$, and
$$\aligned\dd_x(1/f)&=-3x^2/f^2\in F_1\1\OO_{X^{\rm an},0}^{\sD},\\ \dd_x^2(1/f)&=-6x/f^2+18x^4/f^3\\&=\bl(12x^4-6x(y^3+z^3)\br)/f^3\in F_2\1\OO_{X^{\rm an},0}^{\sD}.\endaligned$$
If the equality (2.4.2) holds for $p=2$, then we would get by using the above calculation
$$x(y^3+z^3)\in \Vt^3\!A,$$
since $x^4=f_x^2/9\in \Vt^3\!A$. However, this contradicts (2.2.4), since $x(y^3+z^3)\notin E^2\!A$, where $E$ is defined just before (2.3.2).
\msn
{\bf 2.5.~Calculation of $\alt_f$ in the semi-weighted-homogenous case.} Let $f$ be a weighted homogeneous polynomial of weights $(w_i)$ having an isolated singularity at 0, see (2.2). Then it is well known (see for instance \cite[Chapter 1 or 4]{Sat}) that
$$\alt_f=\msum_i\,w_i.
\leqno(2.5.1)$$
We have the same with $f$ replaced by a semi-weighted-homogeneous polynomial with isolated singularity in the strong sense. The latter means that $f$ is a (finite) linear combination of $f_{\al}$ ($\al\ges 1$) such that each $f_{\al}$ is a linear combination of $\prod_ix_i^{m_i}$ with $\sum_im_iw_i=\al$, and $f_1$ has an isolated singularity at the origin. (Here we can use the finite determinacy of holomorphic functions with isolated singularities to get polynomials.) If there is a semi-weighted-homogeneous polynomial $f=\sum_{\al\ges 1}f_{\al}$ as above, then we have a $\mu$-constant one-parameter family $f^{(\la)}:=\sum_{\al\ges 1}\la^{\al-1}f_{\al}$ for $\la\in\C$ such that $f^{(0)}=f_1$ by replacing $\C$ with a ramified cover to define the $\la^{\al-1}$. (Note that a converse holds by [Va3] using \cite{KaSc,Tj}.)
\sk
The invariance of $\alt_f$ by a $\mu$-constant deformation follows for instance by combining \cite{Br,Ma1,ScSt,Va1,Va2}. More precisely, $\alt_f$ coincides with the minimal spectral number of $f$ in the sense of \cite{St}, since both can be described by using the Brieskorn lattice \cite{Br}, see \cite{Ma1} and \cite{ScSt,Va1}. Moreover the latter number does not change by $\mu$-constant deformations of holomorphic functions with isolated singularities, see \cite{Va2}. (We can also use \cite{DMST} together with \cite{mhm}, see also \cite[Proposition 3.2]{mic}.) This seems to be closely related to \cite[Theorem D]{MP1} where $w_i=1/d$ for any $i$, see also \cite[Theorem 0.9]{mosc}. (It seems that microlocal $V$-filtration is very closely related to the Hodge ideals.)
\bs\bs
\vbox{\centerline{\bf Appendix}
\bsn
We give here some remarks related to papers of Mustata and Popa \cite{MP1,MP2}.}
\msn
{\bf A.1.~A simple proof of the Restriction Theorem.} Let $D$ be a reduced divisor on a smooth complex variety $X$, and $X'$ be a smooth subvariety of $X$ with $D':=D\cap X'$ a divisor on $X'$ (with multiplicities in general). We have the {\it inclusion relation} of subsheaves
$$\I^{(D'_{\rm red},p)}\bl(-(p+1)(D'-D'_{\rm red})\br){}\subset\I^{(D,p)}|\!|_{X'}\q\h{in}\,\,\,\OO_{X'}\q(\forall\,p\ges 0),
\leqno{\rm (A.1.1)}$$
where the right-hand side is the {\it image} of $\I^{(D,p)}\into\OO_X\onto\OO_{X'}$, and $D'_{\rm red}$ is the reduced divisor associated to $D'$, see \cite{MP1,MP2}. Note that $\bl(-(p+1)(D'-D'_{\rm red})\br)$ on the left-hand side cannot be eliminated in general (for instance $f=x^{2a-1}z+y^{2a-1}z+x^ay^a+z^{2a}$ ($a\ges 2$) with $X'=\{z=0\}$).
\sk
We show that the assertion (A.1.1) can be proved quite easily by using {\it only} the $V$-filtration of Kashiwara and Malgrange {\it without using an embedded resolution of} $D\cup X'\subset X$ as in \cite{MP2} (where $X'$ is denoted by $H$).
\sk
By definition (A.1.1) is equivalent to the {\it inclusion relation} of subsheaves
$$F_p\OXDp\subset F_p\OXD|\!|_{X'}\q\h{in}\,\,\,\OO_{X'}(*D')\q(\forall\,p\ges 0),
\leqno{\rm (A.1.2)}$$
where the right-hand side is the image of $F_p\OXD\into\OXD\onto\OXDp$.
\sk
We can prove (A.1.2) as follows (compare with \cite{MP2}).
The assertion is {\it local on} $X$, and is reduced to the case ${\rm codim}_XX'=1$. We may thus assume $X'=\{x_1=0\}$ with $x_1$ a local coordinate of $X$. We have a cartesian diagram
$$\begin{array}{cccccccccccc}
X\setminus D&\buildrel{i'}\over\longleftarrow&X'\setminus D'\\
\llap{$\vcenter{\hbox{$\scriptstyle j$}}$}\big\downarrow&&
\big\downarrow\rlap{$\vcenter{\hbox{$\scriptstyle j'$}}$}\\
X&\buildrel{i}\over\longleftarrow &X'\end{array}$$
Let $V_{x_1}$ be the filtration of Kashiwara and Malgrange along $x_1=0$ as in (1.1). We have the canonical isomorphism of filtered $\D_{X'}$-modules
$$i^!(\OXD,F)=(\OXDp,F).
\leqno{\rm (A.1.3)}$$
The left-hand side is the underlying filtered $\D_{X'}$-module of $i^!j_*\Q_{h,X\setminus D}[d_X]$ (up to a shift of the filtration $F$), and can be defined by the mapping cone (or cokernel in this case) of
$$\Gr_{V_{x_1}}^0x_1:\Gr_{V_{x_1}}^0(\OXD,F)\to\Gr_{V_{x_1}}^1(\OXD,F),
\leqno{\rm (A.1.4)}$$
which underlies the morphism of mixed Hodge modules (up to a shift of the filtration $F$)
$${\rm Var}:\varphi_{x_1,1}\1 j_*\Q_{h,X\setminus D}[d_X-1]\to \psi_{x_1,1}\1 j_*\Q_{h,X\setminus D}(-1)[d_X-1],
\leqno{\rm (A.1.5)}$$
since the mapping cone of the latter {\it represents} the functor $i^!$ (up to a shift of complex) by using a variant of Beilinson's functor $\xi_g$, see \cite[Corollary 2.24]{mhm}. (Here we use \cite{def} as the definition of mixed Hodge modules, since the definition in \cite{mhm} is too complicated.) Then (A.1.3) follows, for instance, from the functorial isomorphism (see \cite[4.4.3]{mhm}):
$$i^!\ssc j_*=j'_*\ssc i^{\prime\,!}.
\leqno{\rm (A.1.6)}$$
Here one can also use the {\it uniqueness} of the open direct image $j'_*$ in the category of mixed Hodge modules (that is, the open direct image is represented by any mixed Hodge module whose underlying $\Q$-complex is the open direct image, see \cite[Proposition 2.11]{mhm}). In this case we need (A.1.6) only for the underlying $\Q$-complexes, and this is easy to verify.
\sk
Note that (A.1.3) implies the isomorphism of $\D_{X'}$-modules
$$\OXDp={\bf L}\1i^*\OXD\,\bl(=C\bl(x_1:\OXD\to\OXD\br)\br),
\leqno{\rm (A.1.7)}$$
since ${\bf L}\1i^*$ for $\D$-modules (that is, for $\OO$-modules with integrable connection) corresponds to $i^*$ under the functor Sol, which is the dual of DR, up to a shift of complex. This means that it corresponds to $i^!$ under the functor DR up to a shift of complex. (We can prove (A.1.7) directly, since $\OXD$ is locally identified with the inductive limit of the inductive system $\{\F_i\}_{i\in\N}$ with $\F_i=\OO_X$ and the transition morphisms $\F_i\to\F_{i+1}$ given by the multiplication by a local defining function $f$ of $D$, and similarly for $\OXDp$.)
\sk
The assertion (A.1.2) is now proved by using the commutative diagram
$$\begin{array}{cccccccccccc}
F_p\Gr_{V_{x_1}}^0\OXD&\buildrel{\gamma_1}\over\longrightarrow&F_p\Gr_{V_{x_1}}^1\OXD\\
\uparrow&&\uparrow\\
F_pV_{x_1}^0\OXD&\buildrel{\gamma_2}\over\longrightarrow&F_pV_{x_1}^1\OXD\\
\downarrow&&\downarrow\\
F_p\OXD&\buildrel{\gamma_3}\over\longrightarrow &F_p\OXD\\
\downarrow&&\downarrow\\
\,\,\OXD&\buildrel{\gamma_4}\over\longrightarrow &\,\,\OXD\end{array}
\leqno{\rm (A.1.8)}$$
Here the $\gamma_i$ are all induced by $x_1:\OXD\to\OXD$, and we consider the {\it images} of ${\rm Coker}\,\gamma_2$, ${\rm Coker}\,\gamma_3$ in ${\rm Coker}\,\gamma_4$. Note that the most upper two vertical morphisms induce a {\it quasi-isomorphism} of mapping cones
$$C(\gamma_2)\simto C(\gamma_1)\q(\h{hence}\,\,\,\,{\rm Coker}\,\gamma_2\simto{\rm Coker}\,\gamma_1),
\leqno{\rm (A.1.9)}$$
by using the {\it acyclicity} of the mapping cone
$$C\bl(x_1:F_pV_{x_1}^{>0}\OXD\simto F_pV_{x_1}^{>1}\OXD\br),
\leqno{\rm (A.1.10)}$$
which follows from \cite[3.2.1.2]{mhp} (where $V_{\al}=V^{-\al}$). In order to verify that we really get the {\it desired\1} morphism from ${\rm Coker}\,\gamma_1$ to ${\rm Coker}\,\gamma_4$ by the above argument, we also consider a similar diagram with $\OXD$ replaced by $\OO_X$ (where the vertical morphisms induce isomorphisms between the ${\rm Coker}\,\gamma_i$, which are all isomorphic to $\OO_{X'}$) together with a canonical morphism from this diagram to the above diagram. (Note that a {\it canonical\1} morphism is not necessarily always the {\it desired\1} one, and a proof is required sometimes. In this case it is enough to show that the morphism is an isomorphism of holonomic $\D$-modules on the complement of $D$. Indeed, the ambiguity of the isomorphism is given by a non-zero constant multiplication in this case.)
\msn
{\bf A.2.~Non-characteristic case.} In the notation of (A.1), the inclusions in (A.1.1--2) become {\it equalities} if $X'$ is {\it non-characteristic} to the $\D_X$-module $\OXD$ (for instance, in the case $X'$ is transversal to any stratum of a Whitney stratification of $D$). In the codimension 1 case, this non-characteristic condition means the following mutually equivalent conditions:
$$(a)\,\,\varphi_{x_1}{\mathbf R}\1 j_*\Q_{X\setminus D}=0;\q(b)\,\,\Gr_{V_{x_1}}^{\al}\OXD=0\,\,\,(\forall\,\al\notin\Z_{>0});\q(c)\,\,\OXD=V_{x_1}^1\OXD.$$
In this case, the morphisms between the cokernels of $\gamma_i$ ($i=1,2,3$) induced by the vertical morphisms of (A.1.8) are all isomorphisms, and these cokernels are canonically isomorphic to $F_p\OXDp$ and to $F_p\OXD|\!|_{X'}$.
\sk
In the non-characteristic case of codimension 1, the above arguments imply that the $F_p\OXD$ are all {\it $x_1$-torsion-free}, since we have ${\rm Ker}\,\gamma_i=0$ for $i=1,\dots,4$. Moreover
$$(F_p/F_q)\OXD:=F_p\OXD/F_q\OXD\,\,\,\h{is also $x_1$-torsion-free for $q<p$},
\leqno{\rm (A.2.1)}$$
where $p$ may be $+\infty$ (that is, $F_p\OXD=\OXD$) by passing to the inductive limit. Indeed, using the canonical isomorphisms
$$V_{x_1}^{\al}(F_p/F_q)\OXD=F_pV_{x_1}^{\al}\OXD/F_qV_{x_1}^{\al}\OXD\q(\al\in\Q),$$
together with \cite[3.2.1.2]{mhp} and condition (c) above, we get the bijectivity of
$$x_1:V_{x_1}^1(F_p/F_q)\OXD\bl(=(F_p/F_q)\OXD\br){}\simto V_{x_1}^2(F_p/F_q)\OXD\,\bl(\subset(F_p/F_q)\OXD\br).$$
\sk
We have local inclusions of $\OO_X$-modules (choosing a local generator $f$)
$$\I^{(D,q)}/\I^{(D,p)}\subset\OO_X/\I^{(D,p)}\cong P_p\OXD/F_p\OXD\subset\OXD/F_p\OXD.$$
So the assertion (A.2.1) implies that
$$\I^{(D,q)}/\I^{(D,p)}\,\,\,\h{are $x_1$-torsion-free for $q<p$.}
\leqno{\rm (A.2.2)}$$
\sk
In the higher codimensional case, let $x_1,\dots,x_r$ be coordinates of $X$ defining the subvariety. Then they form a {\it regular sequence} for $F_p\OXD/F_q\OXD$ and for $\I^{(D,q)}/\I^{(D,p)}$ by an inductive argument, see also \cite{DMST}.
\msn
{\bf Remark.} In the case $0\in D$ is an isolated singularity of $D$, (A.2.2) implies that
$$\dim\bl(\I^{(D,q)}_0/\I^{(D,p)}_0\br)\,\,\,\h{is invariant by $\mu$-constant deformations for $q<p$.}
\leqno{\rm (A.2.3)}$$
\msn
{\bf A.3.~A remark related to the proof of the Subadditivity Theorem.} In order to deduce the Subadditivity Theorem for Hodge ideals from the Restriction Theorem explained in (A.1), we need the isomorphism of mixed Hodge modules
$$j_*(\Q_{h,U}[n])=(j_1)_*\Q_{h,U_1}[n_1]\boxtimes(j_2)_*(\Q_{h,U_2}[n_1]),
\leqno{\rm (A.3.1)}$$
where $j_a:U_a=X_a\setminus D_a\into X_a$ with $D_a$ a reduced divisor on a smooth variety $X_a$ of dimension $n_a$ ($a=1,2$), and $j:U=U_1\times U_2\into X_1\times X_2$. Its proof essentially follows from \cite[3.8.5]{mhm} (asserting the compatibility of affine open direct images with external products) as is noted in \cite{MP2}. Here we have to apply the pull-back functor $\delta^!$ to (A.3.1) under the diagonal morphism $\delta$ in the case $X_1=X_2$ in order to prove the Subadditivity Theorem. The arguments about the {\it stability} of mixed Hodge modules (in the strong sense) by affine open direct images and by external products are not quite complete in \cite{mhm}. Strictly speaking, we should use the arguments as in \cite{def} (using Beilinson's functor \cite{Be} together with the theory of admissible variations of Hodge structure \cite{Ka3}) for the proof of (A.3.1). However, these are not necessary in this special case where we are considering only affine open direct images of Hodge modules with {\it constant coefficients}, and \cite[3.8.5]{mhm} is sufficient. The proof can be reduced essentially to the calculations in the normal crossing case as in \cite[Propositions 3.25--26]{mhm} by using an embedded resolution of the union of the diagonal and the pull-backs of the divisors $D_1,D_2$ by projections.
\msn
{\bf A.4.~Remarks related to the proof of \cite[Proposition 2.4]{MP1}.} It is shown there that we have the following canonical isomorphism in the derived category of $f^{-1}\OO_X$-modules:
$$\omega_Y(*E)\simto\omega_Y(*E)\otimes^{\bf L}_{\D_Y}\D_{Y\to X},
\leqno{\rm (A.4.1)}$$
if $f:Y\to X$ is a proper morphism of smooth complex varieties inducing an isomorphism over the complement of a reduced visor $D$ on $X$, and $E:=(f^*D)_{\rm red}$. The arguments given there seem to be rather too complicated. It can be argued as follows.
\sk
First the assertion must be reduced to the isomorphism in (A.4.1) with $\omega_Y$ replaced by $\D_Y$, that is, to the canonical isomorphism in \cite[Lemma 2.6]{MP1}:
$$\D_Y(*E)\simto\D_Y(*E)\otimes_{\D_Y}\D_{Y\to X},
\leqno{\rm (A.4.2)}$$
where
$$\aligned\D_Y(*E)&:=j'_*\D_{Y\setminus E}\cong\OO_Y(*E)\otimes_{\OO_Y}\D_Y,\\\D_X(*D)&:=j_*\D_{X\setminus D}\cong\OO_X(*D)\otimes_{\OO_X}\D_X,\endaligned
\leqno{\rm (A.4.3)}$$
with $j:X\setminus D\into X$, $j':Y\setminus E\into Y$ natural inclusions, see Remark~(i) below. (Here one can use either the left or right $\D_X$-module structure of $\D_X$ as one wants for the above tensor product, and similarly for $\D_Y$.) The assertion (A.4.1) is local, since there is a canonical morphism. So the reduction can be done by taking a (standard) locally free resolution of $\omega_Y$ over $\D_Y$ to calculate the derived tensor product $\otimes^{\bf L}_{\D_Y}$, and using the {\it filtration} $\sigma$ (see \cite[1.4.7]{De2}) on the locally free resolution (such that its graded quotients are the components of the complex shifted by the degrees) together with Remark~(ii) below. Indeed, (A.4.2) implies that the mapping cone of the morphism representing (A.4.1) is {\it filtered acyclic} for the filtration $\sigma$.
\sk
Note that $\D_Y(*E)$ is {\it flat} over $\D_Y$ (with respect to either the left or right $\D_Y$-module structure of $\D_Y(*E)$ as one wants), since the tensor product with $\D_Y(*E)$ over $\D_Y$ can be identified with the localization along $E$, and is an {\it exact} functor. So ${\bf L}$ over $\otimes_{\D_Y}$ on the right-hand side of (A.4.2) can be omitted. (More precisely, it stays invariant by putting ${\bf L}$ as in \cite[Lemma 2.6]{MP1}, although one does not put it usually when the flatness is known. This is actually related to Remark~(iii) below.)
\sk
In the case of (A.4.2), the morphism is {\it also} $\OO_Y$-linear (using the left $\D_Y$-module structure of $\D_{Y\to X}:=\OO_Y\otimes_{f^{-1}\D_X}f^{-1}\D_X$), and the assertion easily follows by using Remark~(i) below. It is also easy to show the canonical isomorphisms
$$\aligned{\bf R}f_*\D_Y(*E)=f_*\D_Y(*E)&=\D_X(*D)\,\bl(=j_*\D_{X\setminus D}={\bf R}j_*\D_{X\setminus D}\br),\\ f^*\D_X(*D)&=\D_Y(*E)\,\bl(=j'_*\D_{Y\setminus E}\br),\endaligned
\leqno{\rm (A.4.4)}$$
where $f^*\D_X(*D)$ is defined by using either the left or right $f^{-1}\D_X$-module structure of $f^{-1}\D_X(*D)$ as one wants, and the second isomorphism of (A.4.4) is essentially the same as (A.4.2) when the left $f^{-1}\D_X$-module structure is used (since $E=(f^*D)_{\rm red}$).
\msn
{\bf Remarks.} (i) Let $j:X\setminus D\into X$ be the natural inclusion of the complement of a locally principal divisor $D$ on a complex algebraic variety $X$ in general. For any quasi-coherent sheaf $\F$ on $X$, we have a canonical isomorphism of quasi-coherent $\OO_X$-modules
$$\F(*D)\simto j_*j^*\F\,\bl(={\bf R}j_*j^*\F\br),
\leqno{\rm (A.4.5)}$$
where the left-hand side is the localization by a defining function of $D$ (and ${\bf R}j_*$ is the derived direct image of Zariski sheaves using flasque resolutions).
\ms
(ii) For a sheaf of rings $\A$ on a complex variety or on a topological space more generally (for instance $\A=f^{-1}\OO_X$), a morphism in the derived category of $\A$-modules is an isomorphism if and only if it is a {\it quasi-isomorphism} in the generalized sense (here a morphism in the derived category is called a quasi-isomorphism if it induces isomorphisms of all the cohomology sheaves). This condition is also equivalent to the {\it acyclicity} of the mapping cone of the morphism, and stays invariant by forgetting (partially) the action of $\A$ .
\ms
(iii) In the arguments in \cite[Section 2]{MP1}, it seems rather important to clarify in which category the morphism (A.4.2) belongs (if Remark~(ii) above is unused). It must be in the {\it derived category of left $\D_Y(*E)$- and right $f^{-1}\D_X$-bimodules}. In this case, however, the existence of the derived tensor product (taking an appropriate flat resolution) does not seem to be completely trivial, and some explanation may be desirable. It seems better to define the derived tensor product explicitly by taking a (standard) locally free resolution of $\omega_Y$. It is not completely trivial whether the usual construction of flat resolutions for ringed spaces can be applied to the {\it bimodule} case, and it does not seem very clear {\it in which category\1} we can really get the isomorphisms in the proof of \cite[Lemma 2.6]{MP1} (unless Remark~(ii) above is used).
\ms
(iv) We have the {\it bijectivity} of the morphism in \cite[Lemma 2.3]{MP1}:
$$\omega_Y(*E)\simto\omega_Y(*E)\otimes_{\D_Y}\D_{Y\to X},
\leqno{\rm (A.4.6)}$$
which is stated simply as a ``split injection" there.
Indeed, this bijectivity is a corollary of \cite[Proposition 2.4]{MP1} (that is, (A.4.1) above), since the {\it right exactness\1} of tensor products implies
$$\Hc^0\bl(\omega_Y(*E)\otimes^{\bf L}_{\D_Y}\D_{Y\to X}\br){}=\omega_Y(*E)\otimes_{\D_Y}\D_{Y\to X}.
\leqno{\rm (A.4.7)}$$
(It seems rather strange that there is no mention of these assertions in \cite{MP1}.) Note that (A.4.1) is equivalent to (A.4.6) together with the vanishings:
$${\mathcal T}\!or_{\D_Y}^i(\omega_Y(*E),\D_{Y\to X})=0\q(\forall\,i>0).
\leqno{\rm (A.4.8)}$$
\sk
It is also possible to verify (A.4.6) directly. Indeed, shrinking $Y$ and trivializing $\omega_Y$, we have a surjection of right $\D_Y$-modules $\D_Y\to\omega_Y$. This gives the commutative diagram
$$\begin{array}{cccccccccccc}
\D_Y(*E)&\to&\D_Y(*E)\otimes_{\D_Y}\D_{Y\to X}\\ \downarrow&&\downarrow\,\,\,\\\,\omega_Y(*E)&\to&\,\omega_Y(*E)\otimes_{\D_Y}\D_{Y\to X}\end{array}
\leqno{\rm (A.4.9)}$$
where the vertical morphisms are surjective by the right exactness of tensor products, and the top row is bijective by (A.4.2). So the surjectivity of the bottom row follows. Its injectivity can be shown easily by using $j'_*j'{}^*$ as in the proof of \cite[Lemma 2.3]{MP1}.
\msn
{\bf A.5.~Improvement of a vanishing theorem.} We have the vanishing theorem
$$R^{\1 q}\!f_*\Omega_Y^p(\log E)=0\q\h{if}\,\,\,p+q>\dim X,
\leqno{\rm (A.5.1)}$$
under the assumption that $f:Y\to X$ is a proper morphism of complex algebraic varieties, $D$ is a Cartier divisor on $X$, $E:=f^{-1}(D)$ is a divisor with normal crossings on a smooth variety $Y$, and moreover the morphism $Y\setminus E\to X\setminus D$ induced by $f$ is {\it semi-small\,} in the sense of de Cataldo and Migliorini \cite{dCMi}, that is,
$$\aligned&\dim\,(X\setminus D)^k+2k\les\dim Y\q(\forall\,k\ges 0),\\
\h{with}\q\q&(X\setminus D)^k:=\{x\in X\setminus D\mid\dim f^{-1}(x)=k\}.\endaligned$$
Indeed, these conditions imply that the direct image of $(\OO_Y(*E),F)$ as filtered $\D$-module is isomorphic to a filtered $\D_{X''}$-{\it module} with $X''$ a smooth variety containing $X$ locally. (Here we may assume $X''=X$ replacing $X$ locally, since $f$ is not assumed to be surjective.)
\sk
The point is the commutativity of the functor $\DR^{-1}$ with the direct image functor in the filtered derived categories (see \cite[Section 2.3.7]{mhp}) together with the equivalence of categories as in \cite[Proposition 2.2.10]{mhp} where we have the functorial isomorphisms in the derived categories of filtered complexes:
$$\DR\ssc\DR^{-1}\cong id,\q\q\DR^{-1}\ssc\DR\cong id.
\leqno{\rm (A.5.2)}$$
(Note that $\DR$ is denoted by $\widetilde{\rm DR}$ in \cite{mhp} to distinguish it with the de Rham functor to the derived categories of {\it $\C$-complexes}.) We use an $F$-filtered quasi-isomorphism as in \cite[Proposition 3.11(ii)]{mhm} to show the isomorphism in the bounded derived category of filtered right $\D_Y$-modules:
$$\DR^{-1}\bl(\Omega_Y^{\ssb}(\log E),F\br)[d_Y]=\bl(\omega_Y(*E),F\br),
\leqno{\rm (A.5.3)}$$
where the left-hand side is the filtered complex of right $\D_Y$-modules with $i$\1 th component
$$\Omega_Y^{d_Y+i}(\log E)\otimes_{\OO_Y}(\D_Y,F[i])\q(i\in\Z).
\leqno{\rm (A.5.4)}$$
In this {\it constant coefficient} case, we can also use the second isomorphism of (A.5.2) together with the canonical filtered isomorphism (see \cite[Proposition II.3.13(ii)]{De1}):
$$\bl(\Omega_Y^{\ssb}(\log E),F\br)[d_Y]\simto\DR\bl(\omega_Y(*E),F\br).$$
Note that the assertions are shown for the underlying analytic sheaves, and we have to use GAGA, see also Remark~(i) below.
\sk
The above argument implies that the direct image of $\bl(\omega_Y(*E),F\br)$ as filtered right $\D$-module is given by the sheaf-theoretic direct image of
$$\DR^{-1}\bl(\Omega_Y^{\ssb}(\log E),F\br)[d_Y]\otimes_{\D_Y}f^*(\D_X,F),
\leqno{\rm (A.5.5)}$$
whose $i$\1 th component is
$$\Omega_Y^{d_Y+i}(\log E)\otimes_{f^{-1}\OO_X}f^{-1}(\D_X,F[i]).
\leqno{\rm (A.5.6)}$$
This coincides with the filtered complex used in \cite{MP1}. By an argument as in the proof of \cite[Lemma 2.3.6]{mhp}, the sheaf-theoretic direct image of (A.5.5) is identified with $\DR^{-1}$ of the direct image of $(\Omega_Y^{\ssb}(\log E),F\br)[d_Y]$ as filtered differential complex (which is defined by the sheaf-theoretic direct image taking care of the differential appropriately, see \cite[Lemma 2.3.6 and Section 2.3.7]{mhp}).
\msn
{\bf Remarks.} (i) In the simple normal crossing case, we have an \'etale morphism
$$\rho:(Y,0)\to(\C^n,0),$$
induced by $y_i\in\OO_{Y,0}$ ($i\in[1,n]$) which define the irreducible components of $E$ passing through $0\in Y$ for $i\in[1,r]$ so that
$$E=\rho^{-1}E'\q\h{with}\q E':=\{y_1\cdots y_r=0\}\subset\C^n,$$
shrinking $Y$ if necessary, where the $y_i$ are identified with the coordinates of $\C^n$ via $\rho$.
\sk
We have quite explicit descriptions of Hodge ideal and the Hodge and pole order filtrations (using essentially Taylor expansions) on $(\C^n,E')$, but {\it not on} $(Y,E)$ Zariski-locally, since $\OO_{Y,0}\not\cong\OO_{\C^n,0}$ without passing to the {\it Henselization} (by considering their quotient fields). Here we have to use the {\it compatibility} of Hodge ideal (and the Hodge and pole order filtrations) with the {\it pull-backs via \'etale morphisms} (which gives generators of the ideals on $Y$). It may be rather difficult to say that the above description holds ``\'etale locally", since an ``\'etale neighborhood" is usually defined by an \'etale morphism {\it to} $Y$, and not {\it from} $Y$, see \cite{Mi}.
\sk
By a similar reason, it seems safer to give a proof of \cite[Proposition 3.1]{MP1} first on $\C^n$, and then take the pull-back to $Y$ via the \'etale morphism (unless GAGA is used), since it does not seem very clear whether ``Laurent monomials", for instance, work very well on $Y$.
\ms
(ii) We need the uniqueness of the open direct images as in \cite[Proposition 2.11]{mhm} to show the coincidence of the definition of Hodge ideals using an embedded resolution of the divisor as in \cite{MP1} with the one using the Hodge filtration on the open direct image of the structure sheaf as a mixed Hodge module as in the introduction of this paper. Indeed, it is needed to show the isomorphism $(f\ssc j')_*=f_*\ssc j'_*$ for direct image functors in the derived categories of mixed Hodge modules, where $j'$ is as in (A.4.3).

\end{document}